# Asymptotic of the generalized Li's sums which non-negativity is equivalent to the Riemann Hypothesis

**S. K. Sekatskii** (LPMV, Ecole Polytechnique Fédérale de Lausanne, Switzerland)

Recently we have established that the Riemann hypothesis (RH) is equivalent to the non-negativity of "generalized Li's sums" $k_{n,b} := \sum_{\rho}(1-\left(\frac{\rho+b}{\rho-b-1}\right)^n)$ taken over all non-trivial Riemann function zeroes for any real $b$ not equal to -½, see *Ukrainian Math. J.*, **66**, 371-383, 2014; *arXiv*:1304.7895. (Famous Li's criterion corresponds to the case $b=0$ (or $b=1$) here). This makes timely detailed studies of these sums, and in particular also the study of their asymptotic for large $n$. This question, *assuming the truth of RH*, is answered in the present Note. We show that on RH, for large enough $n$, for any real $b \neq -1/2$ one has:

$$k_{n,b} = \sum_{\rho}(1-\left(\frac{\rho+b}{\rho-1-b}\right)^n) = \frac{|2b+1|}{2}n\ln n + \frac{|2b+1|}{2}n(\gamma - 1 - \ln\frac{2\pi}{|2b+1|}) + o(n),$$

where $\gamma$ is Euler – Mascheroni constant.

**Introduction.**

Recently, in Ref. [1] we have established the generalized Bombieri – Lagarias' theorem (see [2] for their original theorem) and the generalized Li's criterion of the truth of the Riemann hypothesis concerning the location of non-trivial zeroes of the Riemann zeta-function (see [3] for the original Li's criterion, and see e.g. [4] for standard definitions and discussion of the general properties of the Riemann zeta-function):

**Theorem 1. (Generalized Bombieri – Lagarias' theorem).** *Let a and $\sigma$ are arbitrary real numbers, $a < \sigma$, and R be a multiset of complex numbers $\rho$ such that*

(i)   $2\sigma - a \notin R$

(ii)  $\sum_{\rho}(1+|\operatorname{Re}\rho|)/(1+|\rho+a-2\sigma|^2) < +\infty$

*Then the following conditions are equivalent*

(a) $\operatorname{Re}\rho \leq \sigma$ *for every* $\rho$;

(b) $\sum_{\rho} \operatorname{Re}(1 - \left(\dfrac{\rho - a}{\rho - 2\sigma + a}\right)^n) \geq 0$ *for n=1, 2, 3...*

(c) *For every fixed $\varepsilon > 0$ there is a positive constant $c(\varepsilon)$ such that*

$$\sum_{\rho} \operatorname{Re}(1 - \left(\dfrac{\rho - a}{\rho - 2\sigma + a}\right)^n) \geq -c(\varepsilon)e^{\varepsilon n}, \ n=1, 2, 3...$$

*If at the same conditions $a > \sigma$ is taken, the point (a) is to be changed to*

(a') $\operatorname{Re}\rho \geq \sigma$ *for every $\rho$,*

*points (b), (c) remain unchanged.*

*If, additionally to the aforementioned conditions, also the following takes place:*

(iii)   *If $\rho \in R$, than $\bar\rho \in R$ with the same multiplicity as $\rho$*

*one can omit the operation of taking the real part in (b), (c), the expressions at question are real. (Here, as usual, $\bar{\rho}$ means a complex conjugate of $\rho$ ).*

**Theorem 2. (Generalized Li's criterion).** *Let a is an arbitrary real number, $a \neq \sigma$, and R be a multiset of complex numbers $\rho$ such that*

(i) $\quad 2\sigma - a \notin R, \ a \notin R$

(ii) $\quad \sum_{\rho}(1+|\operatorname{Re}\rho|)/(1+|\rho+a-2\sigma|^2) < +\infty, \ \sum_{\rho}(1+|\operatorname{Re}\rho|)/(1+|\rho-a|^2) < +\infty$

(iii) *If $\rho \in R$, than $2\sigma - \rho \in R$*

*Then the following conditions are equivalent*

*(a) $\operatorname{Re}\rho = \sigma$ for every $\rho$;*

*(b) $\sum_{\rho} \operatorname{Re}\left(1 - \left(\dfrac{\rho-a}{\rho+a-2\sigma}\right)^n\right) \geq 0$ for any a and n=1, 2, 3...*

*(c) For every fixed $\varepsilon > 0$ and any a there is a positive constant $c(\varepsilon,a)$ such that $\sum_{\rho} \operatorname{Re}\left(1 - \left(\dfrac{\rho-a}{\rho+a-2\sigma}\right)^n\right) \geq -c(\varepsilon,a)e^{\varepsilon n}$, for n=1, 2, 3...*

*If, additionally to the aforementioned conditions, also the following takes place:*

(iv) *If $\rho \in R$, than complex conjugate $\bar{\rho} \in R$ with the same multiplicity as $\rho$*

*one can omit the operation of taking the real part in (b), (c), the expressions at question are real.*

Then, applying the generalized Littlewood theorem about contour integrals of logarithm of an analytical function (see [1, 5-7]), we have established the following equality: for real *a<1/2*

$$\sum_{\rho}\left(1 - \left(\frac{\rho-a}{\rho+a-1}\right)^n\right) = \frac{1}{(n-1)!}\frac{d^n}{dz^n}((z-a)^{n-1}\ln(\xi(z)))\Big|_{z=1-a} \quad (1),$$

and for real *a>1/2*

$$\sum_{\rho}(1-\left(\frac{\rho-a}{\rho+a-1}\right)^n) = -\frac{1}{(n-1)!}\frac{d^n}{dz^n}((z-a)^{n-1}\ln(\xi(z)))|_{z=1-a} \qquad (1a).$$

As follows from Theorem 2, on RH these sums should be non-negative for any real $a$ and any integer $n$, so we have proven the following, see [1]:

**Theorem 3.** *Riemann hypothesis is equivalent to the non-negativity of all derivatives $\frac{1}{(n-1)!}\frac{d^n}{dz^n}((z-a)^{n-1}\ln(\xi(z)))|_{z=1-a}$ for all non-negative integers n and any real a<1/2; correspondingly, it is equivalent also to the non-positivity of all derivatives $\frac{1}{(n-1)!}\frac{d^n}{dz^n}((z-a)^{n-1}\ln(\xi(z)))|_{z=1-a}$ for all non-negative integers n and any real a>1/2.*

Thus to judge the truth of the Riemann hypothesis, certain derivatives of the Riemann xi-function can be estimated at an arbitrary point of the real axis except the point $z=1/2$, not only at the point $z=1$ (or 0) as this was initially formulated by Li [3]. In particular, this can be done far to the right from the point $z=1$, where Riemann zeta-function and its logarithm are defined by absolutely convergent series [4]:

$$\varsigma(z) = \sum_{n=1}^{\infty}\frac{1}{n^z} \qquad (1)$$

$$\ln\varsigma(z) = -\sum_{p}\ln\left(1-\frac{1}{p^z}\right) = \sum_{p}\sum_{n=1}^{\infty}\frac{1}{np^{nz}} = \sum_{n=2}^{\infty}\frac{\Lambda(n)}{\ln n \cdot n^z} \qquad (2)$$

(in (2) we have a sum over primes or use the van Mandgoldt function). This circumstance holds promise to elucidate certain properties of the Riemann function zeroes, and all this make timely the detailed studies of the generalized Li's sums $k_{n,-a} = \sum_{\rho}(1-\left(\frac{\rho-a}{\rho+a-1}\right)^n)$ over non-trivial Riemann function zeroes, and in particular also the study of their asymptotic for large $n$ (which means also the study of asymptotic of corresponding derivatives, see eqs. (1)). This question, *assuming the truth of RH*, is answered in the present Note.

## 2. Asymptotic of the generalized Li's sums and corresponding derivatives assuming the Riemann Hypothesis

Let us now calculate asymptotic of the sums $k_{n,b} = \sum_{\rho}(1-\left(\frac{\rho+b}{\rho-1-b}\right)^n) = \sum_{\rho}(1-\left(\frac{\rho-b-1}{\rho+b}\right)^n)$ over non-trivial Riemann function zeroes for large $n$ (and thus also an asymptotic of equal to them derivatives in (1), $b=-a$) *assuming the Riemann hypothesis*.

**Theorem 4.** *Assume RH. Then for large enough n, for any real* $b \neq -1/2$

$$k_{n,b} = \sum_{\rho}(1-\left(\frac{\rho+b}{\rho-1-b}\right)^n) = \sum_{\rho}(1-\left(\frac{\rho-b-1}{\rho+b}\right)^n) =$$

$$\frac{|2b+1|}{2}n\ln n + \frac{|2b+1|}{2}(\gamma - 1 - \ln\frac{2\pi}{|2b+1|})n + o(n) \qquad (4),$$

where $\gamma$ is Euler – Mascheroni constant.

*Proof.* Proof is a straightforward generalization of the method presented in Coffey paper [8] (see also [9]). Let us first put $b > -1/2$. Using $\rho = 1/2 + iT$, we write for an argument $\vartheta$ of the function $\frac{\rho+b}{\rho-b-1}$:

$$\tan\vartheta = -\frac{(2b+1)T}{T^2 - 1/4 - b - b^2} = -\frac{(2b+1)T}{T^2 - (2b+1)^2/4} \qquad (5).$$

Correspondingly, $\sin\vartheta = -\frac{(2b+1)T}{T^2 + (2b+1)^2/4}$ and $\cos\vartheta = \frac{T^2 - (2b+1)^2/4}{T^2 + (2b+1)^2/4}$; here we used $(T^2 - 1/4 - b - b^2)^2 + (2b+1)^2 T^2 = (T^2 + 1/4 + b + b^2)^2$. Derivative $d\vartheta/dT$ is found from (5): $\frac{d\vartheta}{dT} = \frac{2b+1}{T^2 + (2b+1)^2/4}$, and now we are in a position to calculate the sum at question on RH: $k_{n,b} = \Sigma_\rho (1 - (\frac{\rho+b}{\rho-b-1})^n) = 2\sum_\rho (1 - \cos(n\vartheta_\rho))$ so that, expressed as an integral over the number of non-trivial zeroes $dN$, $k_{n,b} = 2\int_0^\infty (1 - \cos(n\vartheta(T)))dN$. Integrating by parts, we obtain

$$k_{n,b} = 2\int_0^\infty (1 - \cos(n\vartheta(T)))dN = -2n\int_0^\infty \sin(n\vartheta)\frac{d\vartheta}{dT}N(T)dT \qquad (6)$$

and then use the approximations $N(T) = \frac{T}{2\pi}\ln\frac{T}{2\pi} - \frac{T}{2\pi} + O(\ln T)$ [4], $\vartheta = -\frac{2b+1}{T} + O(1/T^3)$, $\frac{d\vartheta}{dT} = \frac{2b+1}{T^2} + O(1/T^4)$ to get

$k_{n,b} = 2n\int_{T_1}^\infty \frac{(2b+1)}{T^2}\sin(\frac{(2b+1)n}{T})N(T)dT + o(n)$ where $T_1 = 14$, say (the first zero lies at $½+i14.1347…$[4]). With the variable change $y = \frac{(2b+1)n}{T}$, we have further

$k_{n,b} = 2\int_0^{(2b+1)n/T_1} \sin(y)N(\frac{(2b+1)n}{y})dy = -\frac{(2b+1)n}{\pi}\int_0^\infty \frac{\sin y}{y}(\ln\frac{2\pi y}{(2b+1)n} + 1)dy + o(n)$, and,

using examples N3.721.1 $\int_0^\infty \frac{\sin y}{y}dy = \frac{\pi}{2}$ and N4.421.1 $\int_0^\infty \ln y \frac{\sin y}{y}dy = -\frac{\pi}{2}\gamma$ from GR book [10], finally obtain

$$k_{n,b} = \sum_\rho (1 - (\frac{\rho+b}{\rho-1-b})^n) = \frac{(2b+1)}{2}n\ln n + \frac{2b+1}{2}(\gamma - 1 - \ln\frac{2\pi}{2b+1})n + o(n).$$

The case $b<-1/2$ is quite similar with changes of signs whenever appropriate, and in this manner we recover the equation (12).

***Remark1.*** Following [8, 9], the sum (derivative) at question can be rewritten, using $\sin n\vartheta = \sin\vartheta \cdot U_{n-1}(\cos\vartheta)$, where $U_k$ is the $k$-th Chebyshev polynomial of the second kind [11], in a rather elegant form

$$k_{n,b} = 2n\int_0^\infty \frac{(2b+1)^2 T}{(T^2+(2b+1)^2/4)^2} U_{n-1}(\frac{T^2-(2b+1)^2/4}{T^2+(2b+1)^2/4}) N(T)dT \qquad (7).$$

We will not use any properties of this polynomial below, but would like to note the next logical step which is the variable change $x = \frac{T^2-(2b+1)^2/4}{T^2+(2b+1)^2/4} = 1 - \frac{(2b+1)^2/2}{T^2+(2b+1)^2/4}$. Clearly, $dx = \frac{T(2b+1)^2}{(T^2+(2b+1)^2/4)^2}dT$ so that

$$k_{n,b} = 2n\int_{-1}^1 U_{n-1}(x)N(x)dx \qquad (8).$$

Using $T = \frac{1}{2}(2b+1)\sqrt{\frac{1+x}{1-x}}$ and limiting ourselves with the $N(T) = \frac{T}{2\pi}\ln\frac{T}{2\pi} - \frac{T}{2\pi} + O(\ln T)$ precision, we may write

$$N(x) = \frac{2b+1}{4\pi}\sqrt{\frac{1+x}{1-x}}\ln(\frac{2b+1}{4\pi}\sqrt{\frac{1+x}{1-x}}) - \frac{2b+1}{4\pi}\sqrt{\frac{1+x}{1-x}} + O(\ln\frac{1+x}{1-x})$$ which is to be

substituted into (8). Note, that integrals of the type $\int_{-1}^1 U_n(x)(1-x)^\alpha (1+x)^\beta dx$ quite naturally appear in some applications of the Chebyshev polynomials; see e. g. example N 7.347.2 of GR book [10].

***Remark 2.*** Note Coffey's suggestion that $o(n)$ terms in the formula (4) for $b=0$ are of the order of $O(n^{1/2+\varepsilon})$ for any $\varepsilon > 0$. Similar property may be suggested for the general case treated here.

***Remark 3.*** It would be interesting to consider an asymptotic, again on RH, of more general sums $k_{n,b,\sigma} = \sum_\rho (1-(\frac{\rho+b}{\rho-2\sigma-b})^n) = \sum_\rho (1-(\frac{\rho-b-1}{\rho+b+2\sigma-1})^n)$ with any real $\sigma \le 1/2$ and any real $b > -\sigma$. According to above Theorem 1, they are certainly positive (because on RH for all non-trivial zeroes $\text{Re}\,\rho \ge \sigma$), but now, if $\sigma \neq 1/2$, the module of any individual summand $|\frac{\rho+b}{\rho-b-2\sigma}|$ is strictly smaller than unity, and earlier used approach with the summing of the terms of the type $2\sum_{\rho}(1-\cos(n\vartheta_\rho))$ cannot be applied any more. We plan to return to this question in the future.